\documentclass[10pt,a4paper,leqno]{amsart}
\usepackage[english]{babel}
\usepackage[T1]{fontenc}
\usepackage{amssymb}
\usepackage{amsfonts}
\usepackage{amsmath}
\usepackage[T1]{fontenc}
\usepackage{eucal}

\begin{document}

%\thispagestyle{empty}
%\twocolumn[
%
%
\null{}
%\vskip3truecm

%\vspace{1ex}

\title{\Large Quasi-nearly subharmonicity and separately quasi-nearly subharmonic functions}
\author{Juhani Riihentaus}
\date{October 16, 2008}
\maketitle
\vspace{-0.35in}
%\begin{multicols}{2}
\begin{center}
{Department of Physics and Mathematics, University of Joensuu\\
P.O. Box 111, FI-80101 Joensuu, Finland \\
juhani.riihentaus@joensuu.fi}
\end{center}
\vspace*{2ex}

%\twocolumn
\noindent{\emph{Abstract:}} Wiegerinck has shown that a separately subharmonic function need not be subharmonic. Improving previous results of Lelong, 
of Avanissian, of Arsove and of us, Armitage and Gardiner gave an almost sharp integrability condition which ensures a  separately subharmonic function 
to be subharmonic. Completing now our recent  counterparts to 
the  cited results of Lelong, Avanissian and Arsove for so called quasi-nearly subharmonic functions, we  present a counterpart to the cited result 
of  Armitage and Gardiner for separately quasi-nearly subharmonic function. This counterpart enables us to slightly improve Armitage's and Gardiner's original 
result, too.  The method we use is a rather straightforward and technical, but still by no means easy, modification of  
Armitage's and Gardiner's argument combined with an 
old argument  of Domar.  

\vspace{0.5ex}


\vspace{3ex}

\noindent{{\textbf{1. Introduction}}}

\vspace{2ex}

\noindent {\textbf{1.1.}} Solving a long standing problem, Wiegerinck [Wi88, Theorem, p.~770], see also Wiegerinck and Zeinstra [WZ91, Theorem~1, p.~246], 
   showed that a separately subharmonic function need not
be subharmonic. On the other hand, Armitage and Gardiner [AG93, Theorem~1, p.~256] showed that a separately subharmonic function $u$ in a domain $\Omega $
in ${\mathbb{R}}^{m+n}$, $m\geq n\geq 2$, is subharmonic provided 
$\phi (\log^+u^+)$ is locally integrable in $\Omega$, where $\phi : [0,+\infty )\rightarrow [0,+\infty )$ is an increasing function such that 
\begin{equation*}\int\limits_1^{+\infty }s^{(n-1)/(m-1)}(\phi (s))^{-1/(m-1)}\,ds<+\infty .\end{equation*}
Armitage's and Gardiner's result includes the previous  results of Lelong [Le45,  Théorème~1~bis, p.~315], of Avanissian 
[Av61, Théorème~9, p.~140], see also [He71, Theorem, p.~31], of  Arsove [Ar66, Theorem~1, p.~622], and of us [Ri89, Theorem~1, p. 69].
Though Armitage's and Gardiner's result is almost sharp, it is, nevertheless, based on Avanissian's result, or, alternatively, on the more 
general results of Arsove and us. See [Ri07$_3$].

In [Ri07$_3$, Proposition~3; Theorem~1, Corollary~1, Corollary~2, Corollary~3; Theorem~2, Corollary; Theorem~3 (Proposition~3.1, p.~57; 
Theorem~3.1, Corollary~3.1, Corollary~3.2, Corollary~3.3, pp.~58-63; Theorem~4.1, Corollary~4.1, pp.~64-65; Theorem~4.2, p.~65)]  we have extended the cited result of Lelong, 
Avanissian, Arsove, and us to the so called quasi-nearly subharmonic functions. The purpose of this paper is to extend also Armitage's and Gardiner's 
result to this more general setup. Our result will, at the same time, give a slight refinement to Armitage's and 
Gardiner's result. The method of proof will be  a rather straightforward and technical, but still by no means easy, modification of Domar's and Armitage's and Gardiner's argument, see 
[Do57, Lemma~1, pp.~431-432 and 430] and [AG93, proof of Proposition~2, \mbox{pp.~257-259,} proof of Theorem~1, pp.~258-259].

\vspace{1ex}

\noindent {\textbf{1.2. Notation.}} Our notation is rather standard, see e.g. [Ri06$_1$], [Ri07$_3$] and [He71]. $m_N$ is the Lebesgue measure 
in the Euclidean space ${\mathbb{R}}^N$, $N\geq 2$. We write $\nu _N$ for the Lebesgue measure of the unit ball $B^N(0,1)$ 
in ${\mathbb{R}}^N$, thus $\nu _N=m_N(B^N(0,1))$. $D$ is a domain of ${\mathbb{R}}^N$.
The complex space ${\mathbb{C}}^n$ is identified with the real space  ${\mathbb{R}}^{2n}$, $n\geq 1$.  
Constants will be denoted by $C$ and $K$. They will be nonnegative and may vary from line to line.

\vspace{3ex}

\noindent{\textbf{2. Quasi-nearly subharmonic functions}} 

\vspace{2ex}

\noindent {\textbf{2.1. Nearly subharmonic functions.}} We recall that an upper semicontinuous function $u:\, D\rightarrow [-\infty ,+\infty )$ is \emph{subharmonic} if 
for all $\overline{B^N(x,r)}\subset D$,
\[u(x)\leq \frac{1}{\nu _N\, r^N}\int\limits_{B^N(x,r)}u(y)\, dm_N(y).\]
The function $u\equiv -\infty $  is considered  subharmonic. 

We say that a function 
$u:\, D\rightarrow [-\infty ,+\infty )$ is \emph{nearly subharmonic}, if $u$ is Lebesgue measurable, $u^+\in {\mathcal{L}}^1_{\textrm{loc}}(D)$, 
and for all $\overline{B^N(x,r)}\subset D$,  
\begin{equation*}u(x)\leq \frac{1}{\nu _N\, r^N}\int\limits_{B^N(x,r)}u(y)\, dm_N(y).\end{equation*}
Observe that in the standard definition of nearly subharmonic functions one uses the slightly stronger  assumption that 
$u\in {\mathcal{L}}^1_{\textrm{loc}}(D)$, see e.g. [He71, p.~14]. However, our above, slightly 
more general definition seems to be  more useful, see  [Ri07$_3$,  Proposition~1~(iii) and Proposition~2~(vi) and (vii) 
(Proposition~2.1~(iii) and Proposition~2.2~(vi), (vii), pp.~54-55)].

\vspace{1ex}

\noindent {\textbf{2.2. Quasi-nearly subharmonic functions.}}  
A Lebesgue measurable function $u:\,D \rightarrow 
[-\infty ,+\infty )$ is \emph{$K$-quasi-nearly subharmonic}, if  $u^+\in{\mathcal{L}}^{1}_{\textrm{loc}}(D)$ and if there is a 
constant $K=K(N,u,D)\geq 1$
such that for all   $\overline{B^N(x,r)}\subset D$,    
\begin{equation} u_M(x)\leq \frac{K}{\nu _N\,r^N}\int\limits_{B^N(x,r)}u_M(y)\, dm_N(y)\end{equation}
for all $M\geq 0$, where $u_M:=\max\{u,-M\}+M$. A function $u:\, D\rightarrow [-\infty ,+\infty )$ is \emph{quasi-nearly subharmonic}, if $u$ is 
$K$-quasi-nearly subharmonic for some $K\geq 1$.

\vspace{1ex}

A Lebesgue  measurable function 
$u:\,D \rightarrow [-\infty ,+\infty )$ is \emph{$K$-quasi-nearly subharmonic n.s. (in the narrow sense)}, if $u^+\in{\mathcal{L}}^{1}_{\textrm{loc}}(D)$ and if there is a 
constant $K=K(N,u,D)\geq 1$
such that for all $\overline{B^N(x,r)}\subset D$,   
\begin{equation} u(x)\leq \frac{K}{\nu _N\,r^N}\int\limits_{B^N(x,r)}u(y)\, dm_N(y).\end{equation}
A function $u:\, D\rightarrow [-\infty ,+\infty )$ is \emph{quasi-nearly subharmonic n.s.}, if $u$ is 
$K$-quasi-nearly subharmonic n.s.  for some $K\geq 1$.

\vspace{1ex}

Quasi-nearly subharmonic functions (perhaps with a different termonology, and sometimes in certain special cases, or just the corresponding generalized mean value inequality (2)) have previously been considered at least in [FS72], [Ku74], [To86], [Ri89],  
[Pa94], [Mi96], [Ri00], [Ri03], [Ri04], [PR08],  [Ri06$_1$], [Ri06$_2$], [Ri07$_3$], [Ko07] and [DP07]. We recall here only that this  function class 
includes, among \mbox{others,}  subharmonic functions, and, more generally,  quasisubharmonic (see e.g. [Le45, p.~309], [Av61, p.~136], [He71, p.~26])
 and also 
nearly subharmonic functions (see e.g. [He71, p.~14]),    also functions satisfying certain natural  growth conditions, especially  
certain eigenfunctions, and  polyharmonic functions. Also, the class of Harnack functions is included, thus, among others, nonnegative harmonic functions 
as well as nonnegative solutions of some elliptic equations. In particular, the partial differential equations associated with quasiregular mappings 
belong to this family of elliptic equations, see Vuorinen [Vu82].  Observe that already Domar  [Do57, p.~430] has pointed out the relevance of the 
class of (nonnegative) quasi-nearly subharmonic functions. For, at least partly,  an even more general function class, see [Do88].

For examples and basic properties of quasi-nearly subharmonic functions, see the above references, especially [PR08] and [Ri07$_3$]. For the sake of 
convenience of the reader we recall the following:  

\begin{itemize}
\item[(i)] \emph{A $K$-quasi-nearly subharmonic function n.s. is $K$-quasi-nearly subharmonic, but not necessarily conversely.}
\item[(ii)] \emph{A nonnegative Lebesgue measurable function is $K$-quasi-nearly subharmonic if and only if it is $K$-quasi-nearly subharmonic n.s.} 
\item[(iii)] \emph{A Lebesgue measurable function is \mbox{$1$-quasi-nearly} subharmonic if and only if it is \mbox{$1$-quasi}-nearly subharmonic n.s. and 
if and only if it is nearly subharmonic (in the sense  defined above).}
\item[(iv)] \emph{If $u:\, D\rightarrow [-\infty ,+\infty )$   is \mbox{$K_1$-quasi-nearly} subharmonic and 
$v:\, D\rightarrow [-\infty ,+\infty )$ is \mbox{$K_2$-quasi-nearly} subharmonic, then $\sup\{u,v\}$ is $\sup\{K_1,K_2\}$-quasi-nearly 
subharmonic in $D$. Especially,  
 $u^+:=\max\{u,0\}$ is $K_1$-quasi-nearly subharmonic in $D$.}
\item[(v)] \emph{Let ${\mathcal{F}}$ be a family of  $K$-quasi-nearly subharmonic (resp. $K$-quasi-nearly subharmonic n.s.) functions in $D$ and let 
$w:=\sup_{u\in {\mathcal{F}}}u$. If $w$ is Lebesgue measurable and} $w^+\in {\mathcal{L}}_{{\textrm{loc}}}^1(D)$, \emph{then $w$ is 
$K$-quasi-nearly subharmonic (resp. $K$-quasi-nearly subharmonic n.s.) in $D$.}
\item[(vi)] \emph{If  $u:\, D\rightarrow [-\infty ,+\infty )$   is quasi-nearly subharmonic n.s., then either $u\equiv -\infty $ or $u$ is finite almost 
everywhere in $D$, and} 
$u\in {\mathcal{L}}^1_{\textrm{loc}}(D)$.
\end{itemize}

\vspace{3ex}

\noindent{{\textbf{3. Lemmas}}} 

\vspace{2ex}

\noindent{\textbf{3.1.}} The following result and its proof is essentially due to Domar [Do57, Lemma~1, pp.~431-432 and 430]. We state the result, however, 
in a more general form, at least seemingly. See also [AG93, p.~258]. 

\vspace{1ex}

\noindent{\textbf{3.2. Lemma.}}  \emph{Let $K\geq 1$. Let $\varphi :[0,+\infty )\rightarrow [0,+\infty )$ and $\psi :[0,+\infty )\rightarrow [0,+\infty )$
 be increasing (strictly or not) functions such that there are $s_0, \,s_1\in {\mathbb{N}}$,  $s_0<s_1$, such that 
\begin{itemize}
\item[{(i)}] the inverse functions $\varphi ^{-1}$ and  $\psi^{-1}$ are defined on $[\inf \{\,\varphi (s_1-s_0),\psi (s_1-s_0)\,\},+\infty )$, 
\item[{(ii)}] $2K(\psi ^{-1}\circ \varphi )(s-s_0)\leq (\psi ^{-1}\circ \varphi )(s)$ for all $s\geq s_1$.
\end{itemize}
Let $u: \,D\rightarrow [0,+\infty )$ be a $K$-quasi-nearly subharmonic function. Suppose that
\begin{equation*}
u(x_j)\geq  (\psi ^{-1}\circ \varphi )(j)
\end{equation*}
for some $x_j\in D$, $j\geq s_1$. If
\begin{equation*}
R_j\geq \left(\frac{2K}{\nu _N}\right)^{1/N}\left[\frac{(\psi ^{-1}\circ \varphi )(j+1)}{(\psi ^{-1}\circ \varphi )(j)}\, m_N(A_j)\right]^{1/N}
\end{equation*}
where
\begin{equation*}
A_j:=\{\, x\in D\, :\, (\psi ^{-1}\circ \varphi )(j-s_0)\leq u(x)<(\psi ^{-1}\circ \varphi )(j+1)\,\},
\end{equation*}
then either $B^N(x_j,R_j)\cap ({\mathbb{R}}^N\setminus D)\ne \emptyset$ or there is $x_{j+1}\in B^N(x_j,R_j)$ such that 
\begin{equation*}
u(x_{j+1})\geq (\psi ^{-1}\circ \varphi )(j+1).
\end{equation*}
}
\vspace{1ex}

\noindent\emph{Proof.} Choose 
\begin{equation*}
R_j\geq \left(\frac{2K}{\nu _N}\right)^{1/N}\left[\frac{(\psi ^{-1}\circ \varphi )(j+1)}{(\psi ^{-1}\circ \varphi )(j)}\, m_N(A_j)\right]^{1/N},
\end{equation*}
and suppose that $B^N(x_j,R_j)\subset D$. Suppose on the contrary that $u(x)<(\psi ^{-1}\circ \varphi )(j+1)$ for all $x\in B^N(x_j,R_j)$. 
Using the assumption (1) (or (2)) we see that 
\begin{align*}
&(\psi  ^{-1}\circ \varphi )(j)\leq u(x_j)\leq \frac{K}{\nu _NR_j^N}\int\limits_{B^N(x_j,R_j)}u(x)\, dm_N(x)\\
&\leq \frac{K}{\nu _NR_j^N}\int\limits_{B^N(x_j,R_j)\cap A_j}u(x)\, dm_N(x)+\frac{K}{\nu _NR_j^N}\int\limits_{B^N(x_j,R_j)\setminus A_j}u(x)\, dm_N(x)
\\
&< \left[\frac{Km_N(B^N(x_j,R_j)\cap A_j)}{\nu _NR_j^N}\frac{(\psi ^{-1}\circ \varphi )(j+1)}{(\psi ^{-1}\circ \varphi )(j)}+
\frac{Km_N(B^N(x_j,R_j)\setminus A_j)}{\nu _NR_j^N}\frac{(\psi ^{-1}\circ \varphi )(j-s_0)}{(\psi ^{-1}\circ \varphi )(j)}\right]
(\psi ^{-1}\circ \varphi )(j)\\
&<(\psi ^{-1}\circ \varphi )(j),
\end{align*}
a contradiction.\hfill \qed

\vspace{1ex}

\noindent{\textbf{3.3.}} The next lemma is a slightly generalized version of Armitage's and Gardiner's result [AG93, Proposition~2, p.~257]. 
The proof of our refinement is -- as already pointed out -- a rather straightforward modification of Armitage's and Gardiner's argument [AG93, proof of Proposition~2, pp.~257-259].

\vspace{1ex}
%\end{document}

\noindent{\textbf{3.4. Lemma.}}  \emph{Let $K\geq 1$.  Let $\varphi :[0,+\infty )\rightarrow [0,+\infty )$ and $\psi :[0,+\infty )\rightarrow [0,+\infty )$
 be increasing functions such that there are $s_0, \,s_1\in {\mathbb{N}}$,  $s_0<s_1$, such  that 
\begin{itemize}
\item[{(i)}] the inverse functions $\varphi ^{-1}$ and  $\psi^{-1}$ are defined on $[\inf \{\,\varphi (s_1-s_0),\psi (s_1-s_0)\,\},+\infty )$, 
\item[{(ii)}] $2K(\psi ^{-1}\circ \varphi )(s-s_0)\leq (\psi ^{-1}\circ \varphi )(s)$ for all $s\geq s_1$,
\item[{(iii)}] $\sum_{j=s_1+1}^{+\infty }\left[\frac{(\psi ^{-1}\circ \varphi )(j+1)}
{(\psi ^{-1}\circ \varphi )(j)}\frac{1}{\varphi (j-s_0)}\right]^{1/(N-1)}<+\infty .$
\end{itemize}
Let $u:\, D\rightarrow [0,+\infty )$ be a $K$-quasi-nearly subharmonic function. Let $\tilde {s}_1\in {\mathbb{N}}$, $\tilde {s}_1\geq s_1$, be arbitrary. 
Then for each $x\in D$ and $r>0$ such that 
$\overline{B^N(x,r)}\subset D$ either
\[u(x)\leq (\psi ^{-1}\circ \varphi )(\tilde {s}_1+1)\]
or
\[\Phi (u(x))\leq \frac{C}{r^N}\int\limits_{B^N(x,r)}\psi (u(y))\, dm_N(y)\]
where $C=C(N,K,s_0)$ and $\Phi :\, [s_2,+\infty )\rightarrow [0,+\infty ),$
\[\Phi (t):=\left(\sum_{j=j_0}^{+\infty }\left[\frac{(\psi ^{-1}\circ \varphi )(j+1)}{(\psi ^{-1}\circ \varphi )(j)}
\frac{1}{\varphi (j-s_0)}\right]^{1/(N-1)}\right)^{1-N},\]
and $j_0\in \{s_1+1,s_1+2,\dots \}$ is such that 
\[(\psi ^{-1}\circ \varphi )(j_0)\leq t<(\psi ^{-1}\circ \varphi )(j_0+1),\]
and $s_2:= (\psi ^{-1}\circ \varphi )(s_1+1)$.
}

\vspace{1ex}

\noindent\emph{Proof.}  Take $x\in D$ and $r>0$ arbitrarily such that $\overline{B^N(x,r)}\subset D$. We may suppose 
that $u(x)>(\psi ^{-1}\circ \varphi )(\tilde {s}_1+1)$. Since $\varphi $
and $\psi $ are increasing, there is an integer $j_0\geq \tilde {s}_1+1$ such that 
\[(\psi ^{-1}\circ \varphi )(j_0)\leq u(x)<(\psi ^{-1}\circ \varphi )(j_0+1).\] 
write $x_{j_0}:=x$, $D_0:=B^N(x_{j_0},r)$ and for each $j\geq j_0$,
\begin{align*}
A_j:=&\{\, y\in D_0\,:\, (\psi ^{-1}\circ \varphi )(j-s_0)\leq u(y)<(\psi ^{-1}\circ \varphi )(j+1)\,\},\\
R_j:=&\left(\frac{2K}{\nu _N}\right)^{1/N}\left[\frac{(\psi ^{-1}\circ \varphi )(j+1)}{(\psi ^{-1}\circ \varphi )(j)}\,
 m_N(A_j)\right]^{1/N}.
\end{align*}
If $B^N(x_{j_0},R_{j_0})\cap ({\mathbb{R}}^N\setminus D_0)\ne \emptyset$, then clearly
\[r<R_{j_0}\leq \sum_{k=j_0}^{+\infty }R_k.\]
On the other hand, if $B^N(x_{j_0},R_{j_0})\subset D_0$, then by Lemma~3.2 there is $x_{j_0+1}\in B^N(x_{j_0},R_{j_0})$ such that 
$u(x_{j_0+1})\geq (\psi ^{-1}\circ \varphi )(j_0+1)$.   Suppose now that for $k=j_0,j_0+1,\dots ,j$,
\[B^N(x_k,R_k)\subset D_0,\,\,   x_{k+1}\in B^N(x_k,R_k)\,\, ({\textrm{this for }} k=j_0, j_0+1, \dots ,j-1),{\textrm{ and }} u(x_k)\geq (\psi ^{-1}\circ \varphi )(k).\]
By Lemma~3.2 there then is $x_{j+1}\in B^N(x_j,R_j)$ such that $u(x_{j+1})\geq (\psi ^{-1}\circ \varphi )(j+1)$. Since $u$ is locally bounded above and 
$(\psi ^{-1}\circ \varphi )(k)\rightarrow +\infty $ as $k\rightarrow +\infty $, we may suppose that 
$B^N(x_{j+1},R_{j+1})\cap ({\mathbb{R}}^N\setminus D_0)\ne \emptyset$. But then 
\begin{displaymath}
r<{\textrm{dist}}(x_{j_0},x_{j_0+1})+{\textrm{dist}}(x_{j_0+1},x_{j_0+2})+\cdots +{\textrm{dist}}(x_{j},x_{j+1})+
{\textrm{dist}}(x_{j+1},{\mathbb{R}}^N\setminus D_0),
\end{displaymath}
thus
\begin{equation}
r< R_{j_0}+R_{j_0+1}+\cdots +R_j+R_{j+1}\leq \sum_{k=j_0}^{+\infty }R_k.
\end{equation}
Using then the  notation 
\[a_k:=\{\, y\in D_0\, :\, (\psi ^{-1}\circ \varphi )(k)\leq u(y)<(\psi ^{-1}\circ \varphi )(k+1)\,\},\]
$k=j_0-s_0, j_0+1-s_0, \dots$, we get from (3):
\begin{align*}
r<&\sum_{k=j_0}^{+\infty }\left(\frac{2K}{\nu _N}\right)^{1/N}\left[\frac{(\psi ^{-1}\circ \varphi )(k+1)}
{(\psi ^{-1}\circ \varphi )(k)} m_N(A_k)\right]^{1/N}\\
<&\left(\frac{2K}{\nu _N}\right)^{1/N}\sum_{k=j_0}^{+\infty }\left(\left[\frac{(\psi ^{-1}\circ \varphi )(k+1)}
{(\psi ^{-1}\circ \varphi )(k)}\frac{1}{\varphi (k-s_0)}\right]^{1/N}\left[\varphi (k-s_0) m_N(A_k)\right]^{1/N}\right)\\
<&\left(\frac{2K}{\nu _N}\right)^{1/N}\left(\sum_{k=j_0}^{+\infty }\left[\frac{(\psi ^{-1}\circ \varphi )(k+1)}
{(\psi ^{-1}\circ \varphi )(k)}\frac{1}{\varphi (k-s_0)}\right]^{1/(N-1)}\right)^{(N-1)/N}
\left[\sum_{k=j_0}^{+\infty }\varphi (k-s_0) m_N(A_k)\right]^{1/N}\\
<&\left(\frac{2K}{\nu _N}\right)^{1/N}\left(\sum_{k=j_0}^{+\infty }\left[\frac{(\psi ^{-1}\circ \varphi )(k+1)}
{(\psi ^{-1}\circ \varphi )(k)}\frac{1}{\varphi (k-s_0)}\right]^{1/(N-1)}\right)^{(N-1)/N}
\left[\sum_{k=j_0}^{+\infty }\int\limits_{A_k}\psi (u(y))dm_N(y)\right]^{1/N}\\
<&\left(\frac{2K}{\nu _N}\right)^{1/N}\left(\sum_{k=j_0}^{+\infty }\left[\frac{(\psi ^{-1}\circ \varphi )(k+1)}
{(\psi ^{-1}\circ \varphi )(k)}\frac{1}{\varphi (k-s_0)}\right]^{1/(N-1)}\right)^{(N-1)/N}
\left(\sum_{k=j_0}^{+\infty }\left[\sum_{j=k-s_0}^{k}\int\limits_{a_j}\psi (u(y))dm_N(y)\right]\right)^{1/N}\\
<&\left[\frac{2(s_0+1)K}{\nu _N}\right]^{1/N}\left(\sum_{k=j_0}^{+\infty }\left[\frac{(\psi ^{-1}\circ \varphi )(k+1)}
{(\psi ^{-1}\circ \varphi )(k)}\frac{1}{\varphi (k-s_0)}\right]^{1/(N-1)}\right)^{(N-1)/N}
\left[\int\limits_{D_0}\psi (u(y))dm_N(y)\right]^{1/N}.
\end{align*}
Thus
\[\Phi (u(x))\leq \frac{C}{r^N}\int\limits_{D_0}\psi (u(y))dm_N(y),\]
where $C=C(N,K,s_0)$ and $\Phi :\, [s_2,+\infty )\rightarrow [0,+\infty ),$
\[\Phi (t):=\left(\sum_{k=j_0}^{+\infty }\left[\frac{(\psi ^{-1}\circ \varphi )(k+1)}{(\psi ^{-1}\circ \varphi )(k)}
\frac{1}{\varphi (k-s_0)}\right]^{1/(N-1)}\right)^{1-N},\]
where $j_0\in \{\, s_1+1,s_1+2, \dots \}$ is such that
\[ (\psi ^{-1}\circ \varphi )(j_0)\leq t<(\psi ^{-1}\circ \varphi )(j_0+1),\]
and $s_2=(\psi ^{-1}\circ \varphi )(s_1+1)$. 

The function $\Phi$ may be extended to the whole interval $[0,+\infty )$, for example as follows:
\begin{displaymath}
\Phi (t):=\begin{cases}\Phi (t), & {\textrm{when }}\, t\geq s_2, \\
 \frac{t}{s_2}\Phi (s_2), & {\textrm{when }}\, 0\leq t<s_2.
\end{cases}\end{displaymath}
\null{}\quad \qed

\vspace{1ex}

\noindent {\textbf{3.5. Remark.}}  Write $s_3:=\sup\{\, s_1+3,(\psi ^{-1}\circ \varphi )(s_1+3)\,\}$, say. (We may suppose that 
$s_3$ is an integer.) 
Suppose that,  in addition to the assumptions (i), (ii), (iii) of Lemma~3.4, also the following assumption  is satisfied:
\begin{itemize}
\item[{(iv)}] \emph{the function 
\[ [s_1+1,+\infty )\ni s\mapsto \frac{(\psi ^{-1}\circ \varphi )(s+1)}
{(\psi ^{-1}\circ \varphi )(s)}\frac{1}{\varphi (t-s_0)}\in {\mathbb{R}}\]
is decreasing.}
\end{itemize}
Then one can replace the function $\Phi \mid [s_3,+\infty )$ by the function $\Phi _1\mid [s_3,+\infty )$, where
\mbox{$\Phi _1(=\Phi _1^{\varphi ,\psi }):\, [0,+\infty )\rightarrow [0,+\infty ),$}
\begin{displaymath}
\Phi_1(t)(=\Phi _1^{\varphi ,\psi }(t)):=\begin{cases}
 \left(\int\limits_{(\varphi  ^{-1}\circ \psi )(t)-2}^{+\infty }\left[\frac{(\psi ^{-1}\circ \varphi )(s+1)}
{(\psi ^{-1}\circ \varphi )(s)}\frac{1}{\varphi (s-s_0)}\right]^{1/(N-1)}ds\right)^{1-N}, & {\textrm{when }}\, t\geq s_3,\\
\frac{t}{s_3}\Phi_1 (s_3), & {\textrm{when }}\, 0\leq t<s_3.
\end{cases}\end{displaymath}

\vspace{1ex}

Similarly, if the function 
\[ [s_1+1,+\infty )\ni s\mapsto \frac{(\psi ^{-1}\circ \varphi )(s+1)}{(\psi ^{-1}\circ \varphi )(s)}\in {\mathbb{R}}\]
is bounded, then in Lemma~3.4 one can replace the function $\Phi \mid [s_3,+\infty )$ by the function $\Phi _2\mid [s_3,+\infty )$, where 
$\Phi _2(=\Phi _2^{\varphi ,\psi }):\, [0,+\infty )\rightarrow 
[0,+\infty ),$
\begin{displaymath}
 \Phi _2(t)(=\Phi _2^{\varphi ,\psi }(t)):=\begin{cases} \left[\int\limits_{(\varphi  ^{-1}\circ \psi )(t)-2}^{+\infty }\frac{ds}{\varphi (s-s_0)^{1/(N-1)}}\right]^{1-N},& 
{\textrm{when }} t\geq s_3, \\
\frac{t}{s_3}\Phi _2(s_3), &{\textrm{when }} 0\leq t<s_3.
\end{cases}
\end{displaymath}
 
\vspace{3ex}

\noindent{{\textbf{4. The condition}}} 

\vspace{2ex}

\noindent{\textbf{4.1.}} Next we propose a counterpart to Armitage's and Gardiner's result [AG93, Theorem~1, p.~256] for quasi-nearly subharmonic 
functions. The proof below follows Armitage's and Gardiner's argument [AG93, proof of Theorem~1, pp.~258-259], but is, at least formally, more general. 
Compare also Corollary~4.7 below. 

\vspace{1ex}

\noindent{\textbf{4.2. Theorem.}} \emph{Let $\Omega $ be a domain in ${\mathbb{R}}^{m+n}$, \mbox{$m\geq n\geq 2$, and let $K\geq 1$.} 
Let $u:\, \Omega \rightarrow [-\infty ,+\infty )$ be  a Lebesgue measurable function. Suppose that the following conditions are satisfied:} 
\begin{itemize}
\item[(a)] \emph{For each $y\in {\mathbb{R}}^n$ the function} 
\[\Omega (y)\ni x\mapsto u(x,y)\in [-\infty ,+\infty )\]
\emph{is $K$-quasi-nearly subharmonic.}
 \item[(b)] \emph{For each $x\in {\mathbb{R}}^m$ the function} 
\[\Omega (x)\ni y\mapsto u(x,y)\in [-\infty ,+\infty )\]
\emph{is $K$-quasi-nearly subharmonic.}
\item[(c)] \emph{There are increasing functions  $\varphi :[0,+\infty )\rightarrow [0,+\infty )$ and $\psi :[0,+\infty )\rightarrow [0,+\infty )$
 and  $s_0, \,s_1\in {\mathbb{N}}$,  $s_0<s_1$, such  that} 
\begin{itemize}
\item[{(c1)}] \emph{the inverse functions $\varphi^{-1}$ and $\psi ^{-1}$ are defined on \mbox{$[\inf \{\, \varphi (s_1-s_0),\psi (s_1-s_0)\,\},
+\infty )$,}}
 \item[{(c2)}] \emph{$2K(\psi ^{-1}\circ \varphi )(s-s_0)\leq (\psi ^{-1}\circ \varphi )(s)$ for all $s\geq s_1$,}
\item[{(c3)}] \emph{the function}
\[ [s_1+1,+\infty )\ni s\mapsto \frac{(\psi ^{-1}\circ \varphi )(s+1)}{(\psi ^{-1}\circ \varphi )(s)}\in {\mathbb{R}}\]
\emph{is bounded,}
\item[{(c4)}] $\int\limits_{s_1}^{+\infty }\frac{s^{(n-1)/(m-1)}}{\varphi (s-s_0)^{1/(m-1)}}ds<+\infty ,$
\item[{(c5)}] $\psi \circ u^+\in {\mathcal{L}}_{\textrm{loc}}^{1}(\Omega )$.
\end{itemize}
\end{itemize}
\emph{Then} $u$ 
\emph{is quasi-nearly subharmonic in $\Omega $}.

\vspace{1ex}

\noindent\emph{Proof.} Recall that  
$s_3=\sup \{\,s_1+3, (\psi ^{-1}\circ \varphi )(s_1+3)\,\}$. Write \mbox{$s_4:=\sup \{\, s_3+s_0, (\varphi ^{-1}\circ \psi  )(s_1+3)\,\}$} and $s_5:=s_4+s_0$,  
say. Clearly, $s_0<s_1<s_2<s_3< s_4< s_5$. (We may suppose that   $s_3$, $s_4$ and $s_5$ are integers.) One may replace $u$ by $u_M:=\sup \{\, u^+,M\,\}$,
where $M=\sup \{\, s_5+3, (\psi ^{-1}\circ \varphi )(s_4+3), (\varphi ^{-1}\circ \psi )(s_4+3)\,\}$, say. We  continue to denote $u_M$ by $u$.

Take $(x_0,y_0)\in \Omega $ and $r>0$ arbitrarily such that $\overline{B^m(x_0,2r)\times B^n(y_0,2r)}\subset \Omega $. 
By [Ri07$_3$, Proposition~3 (Proposition~3.1, p.~57)] (that is, by a counterpart to [Ri89, Theorem~1, p.~69], say)
it is sufficient to show that $u$ is bounded above in $B^m(x_0,r)\times B^n(y_0,r)$.

Take $(\xi ,\eta )\in B^m(x_0,r)\times B^n(y_0,r)$ arbitrarily. In order to apply Lemma~3.4 to the \mbox{$K$-quasi-}nearly subharmonic function $u(\cdot ,\eta )$
in $B^m(\xi ,r)$, check that the assumptions are satisfied. Since (i) and (ii) are clearly satisfied, it remains to show that 
\[\sum_{j=s_1+1}^{+\infty }\left[\frac{(\psi ^{-1}\circ \varphi )(j+1)}{(\psi ^{-1}\circ \varphi )(j)}
\frac{1}{\varphi (j-s_0)}\right]^{1/(m-1)}<+\infty .\]
Because of the assumption~(c3), it is sufficient to show that 
\[\sum_{j=s_1+1}^{+\infty }\frac{1}{\varphi (j-s_0)^{1/(m-1)}}<+\infty .\]
This is seen as follows. Observe first that
\[ \sum_{j=s_1+1}^{+\infty }\frac{1}{\varphi (j-s_0)^{1/(m-1)}}
\leq \int\limits_{s_1}^{+\infty }\frac{ds}{\varphi (s-s_0)^{1/(m-1)}}
\leq \int\limits_{s_1}^{+\infty }\frac{s^{(n-1)/(m-1)}}{\varphi (s-s_0)^{1/(m-1)}}ds<+\infty .\]
We know that $u(\xi ,\eta )\geq s_4$. Therefore it  follows from Lemma~3.4 and Remark~3.5 that
\begin{equation}
\Phi_2 (u(\xi ,\eta ))=\left[\int\limits^{+\infty }_{(\varphi ^{-1}\circ \psi )(u(\xi ,\eta ))-2}\frac{ds}{\varphi (s-s_0)^{1/(m-1)}}\right]^{1-m}
\leq \frac{C}{r^m}\int\limits_{B^m(\xi ,r)}\psi (u(x,\eta ))dm_m(x).
\end{equation}
Recall that here $\Phi_2(=\Phi _2^{\varphi ,\psi }) :\, [0,+\infty )\rightarrow [0,+\infty )$,
\begin{displaymath}
\Phi_2(t)(=\Phi _2^{\varphi ,\psi }(t)):=\begin{cases}
 \left[\int\limits_{(\varphi ^{-1}\circ \psi )(t)-2}^{+\infty }\frac{ds}{\varphi (s-s_0)^{1/(m-1)}}\right]^{1-m}, & {\textrm{when }}\, t\geq s_3,\\
\frac{t}{s_3}\Phi_2 (s_3), & {\textrm{when }}\, 0\leq t<s_3. 
\end{cases}\end{displaymath}

Take then the integral mean values of both sides of (4) over $B^n(\eta ,r)$:
\begin{equation}\begin{split}\frac{C}{r^n}\int\limits_{B^n(\eta ,r)}\Phi _2(u(\xi ,y))dm_n(y)
&\leq \frac{C}{r^n}\int\limits_{B^n(\eta ,r)}[\frac{C}{r^m}\int\limits_{B^m(\xi ,r)}\psi (u(x,y))dm_m(x)]dm_n(y)\\
&\leq \frac{C}{r^{m+n}}\int\limits_{B^m(\xi ,r)\times B^n(\eta ,r)}\psi (u(x,y)) dm_{m+n}(x,y)\\  
&\leq \frac{C}{r^{m+n}}\int\limits_{B^m(x_0 ,2r)\times B^n(y_0 ,2r)}\psi (u(x,y)) dm_{m+n}(x,y).
\end{split}\end{equation}

In order to apply Lemma~3.4 (and Remark~3.5) once more, define  $\psi _1:\, [0,+\infty )\rightarrow [0,+\infty )$, 
$\psi _1(t):=\Phi _2(t)$, where $\Phi_2:\,[0,+\infty )\rightarrow [0,+\infty )$ is as above.
%\begin{displaymath}
%\tilde {\Phi }_2(t):=\begin{cases}\Phi _2(t)=\left[\int\limits_{(\varphi ^{-1}\circ \psi )(t)-2}^{+\infty }\frac{ds}{\varphi (s-s_0)^{1/(m-1)}}\right]^{1-m},
%& {\textrm{when }}\, t\geq s_3, \\
% \frac{t}{s_4}\Phi _2(s_4)=\frac{t}{s_4}\tilde {\Phi }_2(s_),  & {\textrm{when }}\, 0\leq t<s_4.
%\end{cases}\end{displaymath}
Define $\varphi _1:\, [0,+\infty )\rightarrow  [0,+\infty )$,
\begin{displaymath}
\varphi _1(t):=\begin{cases}\frac{t}{s_3}\psi_1 ((\psi ^{-1}\circ \varphi )(s_3))=\frac{t}{s_3}\Phi _2(\psi ^{-1}(\varphi (s_3))), 
& {\textrm{when }}\, 0\leq t<s_3, \\
 \psi _1((\psi ^{-1}\circ \varphi )(t))=\Phi _2(\psi ^{-1}(\varphi (t))), & {\textrm{when }}\, t\geq s_3.
\end{cases}\end{displaymath}
It is straightforward to see that both $\psi _1$ and $\varphi _1$ are strictly increasing and continuous.
Observe also that for $t\geq s_4$, say,  
\begin{equation}\begin{split}
\varphi _1(t)&=\Phi _2((\psi ^{-1}\circ \varphi )(t))
=\left[\int\limits_{(\varphi ^{-1}\circ \psi )((\psi ^{-1}\circ \varphi )(t))-2}^{+\infty }
\frac{ds}{\varphi (s-s_0)^{1/(m-1)}}\right]^{1-m}\\
&=\left[\int\limits_{t-2}^{+\infty }
\frac{ds}{\varphi (s-s_0)^{1/(m-1)}}\right]^{1-m}.\end{split}\end{equation}

Check then that the assumptions of  Lemma~3.4 (and Remark~3.5) are fullfilled. Write $\tilde{s}_0:=s_0$ and  $\tilde{s}_1:=s_4$. 
The assumption (i) is clearly satisfied for $\psi _1$ and $\varphi _1$.  Then  for all $s\geq \tilde {s}_1$,
\[\varphi _1(t)=\psi _1((\psi ^{-1}\circ \varphi )(t))\Leftrightarrow (\psi_1 ^{-1}\circ \varphi _1)(t)=(\psi ^{-1}\circ \varphi )(t).\]
Thus  also the assumption (ii) is  satisfied. It remains to show that  
\[\sum_{j=s_4+1}^{+\infty }\left[\frac{(\psi_1 ^{-1}\circ \varphi_1 )(j+1)}{(\psi_1 ^{-1}\circ \varphi_1 )(j)}
\frac{1}{\varphi_1 (j-s_0)}\right]^{1/(n-1)}<+\infty ,\]
say. It is surely sufficient to show that
\begin{equation}\int\limits_{s_5+s_0+2}^{+\infty }\frac{ds}{{\varphi_1 (s-s_0)}^{1/(n-1)}}<+\infty .\end{equation}

Define $F:\, [s_5,+\infty )\times [s_5+s_0+2,+\infty )\rightarrow [0,+\infty )$,
\begin{displaymath}
F(s,t):=\begin{cases}0, 
& {\textrm{when }}\, s_5\leq s<t-s_0-2, \\
 \varphi (s-s_0)^{-1/(m-1)}, & {\textrm{when }}\, s_5+s_0+2\leq t-s_0-2\leq s.
\end{cases}\end{displaymath}

Suppose that $m>n$ and write $p:=(m-1)/(n-1)$. Using Minkowski's Inequality, see e.g. [Ga07, p.~14], one obtains, with the aid of (6),
\begin{align*}
&\left[\int\limits_{s_5+s_0+2}^{+\infty}\frac{dt}{\varphi _1(t-s_0)^{1/(n-1)}}\right]^{(n-1)/(m-1)}= 
\left[\int\limits_{s_5+s_0+2}^{+\infty }\left(\left[\int\limits_{t-s_0-2}^{+\infty }\frac{ds}{\varphi (s-s_0)^{1/(m-1)}}\right]^{1-m}\right)^{-1/(n-1)}
dt\right]^{(n-1)/(m-1)}\\
&=\left(\int\limits_{s_5+s_0+2}^{+\infty }\left[\int\limits_{t-s_0-2}^{+\infty }\frac{ds}{\varphi (s-s_0)^{1/(m-1)}}\right]^{(m-1)/(n-1)}
dt\right)^{(n-1)/(m-1)}\\
&=\left(\int\limits_{s_5+s_0+2}^{+\infty }\left[\int\limits_{s_5}^{+\infty }F(s,t)ds\right]^{(m-1)/(n-1)}dt\right)^{(n-1)/(m-1)}
\leq \int\limits_{s_5}^{+\infty }\left[\int\limits_{s_5+s_0+2}^{+\infty }F(s,t)^{(m-1)/(n-1)}dt\right]^{(n-1)/(m-1)}ds\\
&\leq \int\limits_{s_5}^{+\infty }\left[\int\limits_{s_5+s_0+2}^{s+s_0+2}\frac{dt}{\varphi (s-s_0)^{1/(n-1)}}\right]^{(n-1)/(m-1)}ds
\leq \int\limits_{s_5}^{+\infty }\frac{[(s+s_0+2)-(s_5+s_0+2)]^{(n-1)/(m-1)}}{\varphi (s-s_0)^{1/(m-1)}}ds\\
&\leq \int\limits_{s_5}^{+\infty }\frac{(s-s_5)^{(n-1)/(m-1)}}{\varphi (s-s_0)^{1/(m-1)}}ds
\leq \int\limits_{s_5}^{+\infty }\frac{s^{(n-1)/(m-1)}}{\varphi (s-s_0)^{1/(m-1)}}ds
<+\infty .
\end{align*}

The case $m=n$ is considered similarly,  just replacing  Minkowski's Inequality with Fubini's Theorem.

Now we can apply Lemma~3.4 (and Remark~3.5) to the left hand side of (5). Recall that $\tilde{s}_0=s_0$,  $\tilde{s}_1=s_4$, 
$\tilde {s}_3:=\sup \{ \, \tilde{s}_1+3, (\psi_1^{-1}\circ \varphi _1 )(\tilde{s}_1+3)\,\}$, 
and \mbox{$\tilde{s}_4:=\sup \{ \, \tilde{s}_3+\tilde{s}_0, (\varphi _1^{-1}\circ \psi_1 )(\tilde{s}_1+3)\,\}$.}  
(Here and below, in the previous definitions just replace the functions
$\varphi $ and $\psi $ with the functions $\varphi _1$ and $\psi _1$, respectively.) Write moreover  
\mbox{$s_4^*:=\sup \{\, \tilde{s}_4,(\psi^{-1}\circ \varphi )(s_4)\,\}$,} say. 
Since $u(\xi ,\eta )\geq M\geq s_4^*\geq \tilde{s}_4$ for all $(\xi ,\eta )\in B^m(x_0,r)\times B^n(y_0,r)$,
we  obtain, using (5), 
\begin{equation}\begin{split}\Psi (u(\xi ,\eta ))&=\left[\int\limits_{(\varphi_1 ^{-1}\circ \psi_1)(u(\xi ,\eta ))-2}^{+\infty }
\frac{ds}{\varphi_1 (s-s_0)^{1/(n-1)}}\right]^{1-n}\\
&\leq \frac{C}{r^n}\int\limits_{B^n(\eta  ,r)}\Phi _2(u(\xi ,y))dm_n(y)\\
&\leq \frac{C}{r^{m+n}}\int\limits_{B^m(x_0,2r)\times B^n(y_0,2r)}\psi (u(x,y))dm_{m+n}(x,y).
\end{split}\end{equation}
Here now  $\Psi (=\Phi _2^{\varphi _1,\psi _1}) :\, [0,+\infty )\rightarrow [0,+\infty )$,
\begin{displaymath}
\Psi (t)(= \Phi _2^{\varphi _1,\psi _1}(t)):=\begin{cases}\left[\int\limits_{(\varphi_1 ^{-1}\circ \psi_1)(t)-2}^{+\infty }\frac{ds}{\varphi_1 (s-s_0)^{1/(n-1)}}\right]^{1-n},
& {\textrm{when }}\, t\geq \tilde {s}_3, \\
 \frac{t}{\tilde {s}_3}\Psi (\tilde {s}_3),  & {\textrm{when }}\, 0\leq t<\tilde {s}_3,
\end{cases}\end{displaymath}
see Remark~3.5 above.
From (8),  from the   facts  that 
$(\varphi_{1} ^{-1}\circ \psi _1)(t)=(\varphi ^{-1}\circ \psi )(t)\rightarrow +\infty $ as $t\rightarrow +\infty $,  from (7), and from 
the fact that
\[\int\limits_{B^m(x_0,2r)\times B^n(y_0,2r)}\psi (u(x,y))dm_{m+n}(x,y)<+\infty ,\] 
one sees that $u$ must be bounded above in $B^m(x_0,r)\times B^n(y_0,r)$, concluding the proof.
\null{}\quad \qed 

\vspace{1ex}

\noindent{\textbf{4.3. Corollary.}} \emph{Let $\Omega $ be a domain in ${\mathbb{R}}^{m+n}$, \mbox{$m\geq n\geq 2$, and let $K\geq 1$.} 
Let $u:\, \Omega \rightarrow [-\infty ,+\infty )$ be  a Lebesgue measurable function. Suppose that the following conditions are satisfied:} 
\begin{itemize}
\item[(a)] \emph{For each $y\in {\mathbb{R}}^n$ the function} 
\[\Omega (y)\ni x\mapsto u(x,y)\in [-\infty ,+\infty )\]
\emph{is $K$-quasi-nearly subharmonic.}
 \item[(b)] \emph{For each $x\in {\mathbb{R}}^m$ the function} 
\[\Omega (x)\ni y\mapsto u(x,y)\in [-\infty ,+\infty )\]
\emph{is $K$-quasi-nearly subharmonic.}
\item[(c)] \emph{There is a strictly  increasing surjection  $\varphi :[0,+\infty )\rightarrow [0,+\infty )$ such that} 
\begin{itemize}
\item[{(c1)}] $\int\limits_{s_0+1}^{+\infty }\frac{s^{(n-1)/(m-1)}}{\varphi (s-s_0)^{1/(m-1)}}ds<+\infty $ \emph{for some $s_0\in {\mathbb{N}}$,} 
\item[{(c2)}] $\varphi (\log^+u^+)\in {\mathcal{L}}_{\textrm{loc}}^{1}(\Omega )$.
\end{itemize}
\end{itemize}
\emph{Then} $u$ 
\emph{is quasi-nearly subharmonic in $\Omega $}.

\vspace{1ex}

\noindent\emph{Proof.} Just choose $\psi =\varphi \circ \log^+$ and apply Theorem~4.2. 
\null{}\quad \qed 

\vspace{1ex}

\noindent {\textbf{4.4. Remark.}} One sees easily that the condition  (c1) can be replaced by the condition
\begin{itemize}
\item[{(c1')}] $\int\limits_{1}^{+\infty }\frac{s^{(n-1)/(m-1)}}{\varphi (s)^{1/(m-1)}}ds<+\infty $. 
\end{itemize}
\vspace{1ex}

\noindent{\textbf{4.5. Corollary.}} \emph{Let $\Omega $ be a domain in ${\mathbb{R}}^{m+n}$, \mbox{$m\geq n\geq 2$, and let $K\geq 1$.} 
Let $u:\, \Omega \rightarrow [-\infty ,+\infty )$ be  a Lebesgue measurable function. Suppose that the following conditions are satisfied:} 
\begin{itemize}
\item[(a)] \emph{For each $y\in {\mathbb{R}}^n$ the function} 
\[\Omega (y)\ni x\mapsto u(x,y)\in [-\infty ,+\infty )\]
\emph{is $K$-quasi-nearly subharmonic.}
 \item[(b)] \emph{For each $x\in {\mathbb{R}}^m$ the function} 
\[\Omega (x)\ni y\mapsto u(x,y)\in [-\infty ,+\infty )\]
\emph{is $K$-quasi-nearly subharmonic.}
\item[(c)] \emph{There is a strictly  increasing surjection  $\varphi :[0,+\infty )\rightarrow [0,+\infty )$ such that} 
\begin{itemize}
\item[{(c1)}] $\int\limits_{s_0+1}^{+\infty }\frac{s^{(n-1)/(m-1)}}{\varphi (s-s_0)^{1/(m-1)}}ds<+\infty $ \emph{for some $s_0\in {\mathbb{N}}$,} 
\item[{(c2)}] $\varphi (\log(1+(u^+)^r))\in {\mathcal{L}}_{\textrm{loc}}^{1}(\Omega )$ \emph{for some} $r>0$.
\end{itemize}
\end{itemize}
\emph{Then} $u$ 
\emph{is quasi-nearly subharmonic in $\Omega $}.

\vspace{1ex}

\noindent\emph{Proof.} It is easy to see that the assumptions of Theorem~4.2 are satisfied. We leave the details to the reader.\null{}\quad \qed 

\vspace{1ex}

\noindent{\textbf{4.6.}} Next our  slight improvement to Armitage's and Gardiner's original result:

\vspace{1ex}
 
\noindent{\textbf{4.7. Corollary.}} \emph{Let $\Omega $ be a domain in ${\mathbb{R}}^{m+n}$, \mbox{$m\geq n\geq 2$.} 
Let $u:\, \Omega \rightarrow [-\infty ,+\infty )$ be such that the following conditions are satisfied:} 
\begin{itemize}
\item[(a)] \emph{For each $y\in {\mathbb{R}}^n$ the function} 
\[\Omega (y)\ni x\mapsto u(x,y)\in [-\infty ,+\infty )\]
\emph{is subharmonic.}
 \item[(b)] \emph{For each $x\in {\mathbb{R}}^m$ the function} 
\[\Omega (x)\ni y\mapsto u(x,y)\in [-\infty ,+\infty )\]
\emph{is  subharmonic.}
\item[(c)] \emph{There is a strictly  increasing surjection  $\varphi :[0,+\infty )\rightarrow [0,+\infty )$ such that} 
\begin{itemize}
\item[{(c1)}] $\int\limits_{1}^{+\infty }\frac{s^{(n-1)/(m-1)}}{\varphi (s)^{1/(m-1)}}ds<+\infty $, 
\item[{(c2)}] $\varphi (\log^+[(u^+)^r])\in {\mathcal{L}}_{\textrm{loc}}^{1}(\Omega )$ \emph{for some} $r>0$.
\end{itemize}
\end{itemize}
\emph{Then} $u$ 
\emph{is subharmonic in $\Omega $}.

\vspace{1ex}

\noindent\emph{Proof.} By [Ri07$_3$, Proposition~2~(v), (vi), and Proposition~1~(iv) (Proposition~2.2~(v), (vi), and Proposition~2.1~(iv), p.~55)], see also [Ri06$_1$, Lemma~2.1, p.~32], $(u^+)^r$ satisfies the assumptions of Corollary~4.3. Thus  $(u^+)^r$ 
is quasi-nearly subharmonic in $\Omega $, and therefore e.g. by [Ri07$_3$, Proposition~2~(iii) (Proposition~2.2~(iii), p.~55)] locally bounded above. Hence also $u$ is locally 
bounded above, and thus 
subharmonic in $\Omega $, by [Ri89, Theorem~1, p.~69], say.
\null{}\quad \qed 

\vspace{1ex}

\noindent{{\textbf{Added 16 Oct 2008:}}} For a further, still slightly improved version of the above Corollary~4.7, see [Ri08, Corollary~3.3.3]. 

\vspace{2ex}

\flushleft{\noindent\textbf{References}}

\vspace{1ex}
\begin{flushleft}
\begin{enumerate}
%\item[{[ABR01]}] Axel, S., Bourdon, P., Ramey, W.  ``Harmonic Function Theory'',  
%Springer-Verlag, New-York, 2001 (Second Edition).
\item[{[AG93]}] Armitage, D.H., Gardiner, S.J.
 ``Conditions for separately subharmonic functions to be subharmonic'',  
Potential Anal.,
{\textbf{2}}, No.~3 (1993), 255--261.
\item[{[AG01]}] Armitage, D.H., Gardiner, S.J.
 ``Classical Potential Theory'',  
Springer-Verlag, London, 2001.
\item[{[Ar66]}] Arsove, M.G. ``On subharmonicity of doubly subharmonic functions'',  
Proc. Amer. Math. Soc.,
{\textbf{17}} (1966), 622--626.
\item[{[Av61]}] Avanissian, V.
 ``Fonctions plurisousharmoniques et fonctions doublement sousharmoniques'',  
Ann. Sci. École  Norm. Sup., {\textbf{78}} (1961), 101--161.
%\item[{[Av67]}] Avanissian, V.
% ``Sur l'harmonicité des fonctions séparément harmoniques'', in: 
%Séminaire de Probabilités (Univ. Strasbourg, Février 1967), {\textbf{1}} (1966/1967), pp.~101--161, Springer, Berlin, 1967.
%\item[{[Br38]}] Brelot, ~M. 
% ``Sur le potentiel et les suites de fonctions sousharmoniques'',  
%C.R. Acad. Sci., {\textbf{207}} (1938), 836--838.
%\item[{[Br69]}] Brelot, M.
 %``Éléments de la Théorie Classique du Potentiel'',  
%Centre de Documentation Universitaire, Paris, 1969 (Third Edition).
%\item[{[CZ61]}] Calderon, A.P., Zygmund, A. ``Local properties of solutions of elliptic partial differential equations'', 
%Studia Math., {\textbf{20}} (1961), 171--225.
%\item[{[CS93]}] Cegrell, U.,  Sadullaev, A. ``Separately subharmonic functions'', Uzbek. Math. J., {\textbf{1}} (1993), 78--83.
%\item[{[Di60]}] Dieudonné, J.
% ``Foundations of Modern Analysis'',  
%Academic Press, New York, 1960.
\item[{[DP07]}] Djordjevi\'c, O., Pavlovi\'c, M. ``${\mathcal{L}}^p$-integrability of the maximal function of a polyharmonic function'', 
J. Math. Anal. Appl., {\textbf{336}}, No.~1 (2007), 411--417.
\item[{[Do57]}] Domar, Y. ``On the existence of a largest subharmonic minorant of a given function'', 
Arkiv f\"or Matematik, {\textbf{3}}, nr. 39 (1957), 429--440.
\item[{[Do88]}] Domar, Y. ``Uniform boundedness in families related to subharmonic functions'', 
J. London Math. Soc. (2), {\textbf{38}}, No.~3 (1988), 485--491.
\item[{[FS72]}] Fefferman, C., Stein, E.M.  ``H$^p$ spaces of several variables'',
Acta Math., \textbf{129} (1972),   \mbox{137--192.}
\item[{[Ga07]}] Garnett,~J.B. ``Bounded Analytic Functions'', 
Springer-Verlag, New York, 2007.
%\item[{[Hel69]}] Helms, L.L. ``Introduction to Potential Theory'', 
%Wiley-Interscience, New York, 1969.
\item[{[He71]}] Hervé, M. ``Analytic and Plurisubharmonic Functions in Finite and Infinite Dimensional Spaces'', 
Lecture Notes in Mathematics,  198, Springer-Verlag, Berlin, 1971.  
%\item[{[Im90]}] Imomkulov, S.A. ``Separately subharmonic functions'' (in Russian),  
%Dokl. USSR, {\textbf{2}} (1990), 8--10.
\item[{[Ko07]}] Koji\'c, V. ``Quasi-nearly subharmonic functions and conformal mappings'', 
Filomat., {\textbf{21}}, No. 2 (2007), 243--249.
\item[{[Ku74]}]Kuran, \"U. ``Subharmonic behavior of $\mid h\mid ^p$, ($p>0$, $h$ harmonic)'',
J. London Math. Soc. (2), {\textbf{8}} (1974),  \mbox{529--538.}
%\item[{[KT96]}] Ko\l odziej, S., Thorbi$\ddot {\textrm{o}}$rnson, J. ``Separately harmonic and subharmonic functions'',  
%Pot. Anal.,
%{\textbf{5}} (1996), 463--466.
\item[{[Le45]}] Lelong, P. ``Les fonctions plurisousharmoniques'',  
Ann. Sci. École Norm. Sup.,
{\textbf{62}} (1945), 301--338.
%\item[{[Le61]}] Lelong, P. ``Fonctions plurisousharmoniques et fonctions analytiques de variables réelles'',  
%Ann. Inst. Fourier, Grenoble, 
%{\textbf{11}} (1961), 515--562.
\item[{[Le69]}] Lelong, P. ``Plurisubharmonic Functions and Positive Differential Forms'',  
Gordon and Breach, London, 1969.
%\item[{[LL01]}] Lieb, E.H., Loss, M. ``Analysis'',
 %Graduate Studies in Mathematics, 14, American Mathematical Society, Providence, Rhode Island, 2001.
%\item[{[Ma95]}] Mattila,~P. ``Geometry of Sets and Measures in Euclidean Spaces'',
%Cambridge studies in advanced mathematics, 44, Cambridge University Press, Cambridge,  1995.
\item[{[Mi96]}] Mizuta,~Y. ``Potential Theory in Euclidean Spaces'',
 Gaguto International Series, Mathematical Sciences and Applications, 6, Gakk$\bar{{\textrm{o}}}$tosho Co., Tokyo, 1996.
\item[{[Pa94]}] Pavlovi\'c, M. ``On subharmonic behavior and oscillation of functions on balls in ${\mathbb{R}}^n$'', 
Publ. Inst. Math. (Beograd), {\textbf{55 (69)}} (1994), 18--22.
\item[{[PR08]}] Pavlovi\'c, M., Riihentaus, J.,  ``Classes of quasi-nearly subharmonic functions'', Potential Anal. {\textbf{29}}, No.~1 (2008), 89--104.
%\item[{[Pl70]}] du Plessis, N. ``An Introduction to Potential Theory'', Oliver \& Boyd, Edinburgh, 1970.
%\item[{[Ra37]}] Radó, T. ``Subharmonic Functions'',
% Springer, Berlin, 1937.
%\item[{[Ri84]}] Riihentaus, J. ``On the extension of separately hyperharmonic functions and H$^p$-functions'',
%Michigan Math. J., {\textbf{31}} (1984),  \mbox{99--112.}
\item[{[Ri89]}] Riihentaus, J. ``On a theorem of Avanissian--Arsove'',
Expo.  Math., {\textbf{7}}, No.~1 (1989),  \mbox{69--72.}
\item[{[Ri00]}] Riihentaus, J.  ``Subharmonic functions: non-tangential and 
tangential boundary
behavior'', in:  Conference Function Spaces, Differential Operators and Nonlinear Analysis (FSDONA'99), Sy\"ote, Pudasj\"arvi, Finland, June 10-16, 1999,  Proceedings,    V.~Mustonen, J.~R\'akosnik (eds.), 
 Math. Inst., Czech Acad. Science,  Praha, 2000, \mbox{pp. 229--238.} \mbox{(ISBN 80-85823-42-X)}
\item[{[Ri01]}] Riihentaus, J.  ``A generalized mean value inequality for subharmonic functions'', Expo. Math., {\textbf{19}}, No.~2 (2001),  \mbox{187-190}.
\item[{[Ri03]}] Riihentaus, J. ``A generalized mean value inequality for subharmonic functions and applications'', arXiv:math.CA/0302261 v1 21 Feb 2003 
(v2 1 Nov 2006). 
\item[{[Ri04]}] Riihentaus, J.  ``Weighted boundary behavior and nonintegrability of subharmonic functions'', in: International Conference on 
Education and Information Systems: Technologies and Applications (EISTA'04), Orlando, Florida, USA, July 21-25, 2004,  Proceedings,
M.~Chang, Y-T.~Hsia, F.~Malpica, M.~Suarez, A.~Tremante, F.~Welsch (eds.), vol.~II, 2004, \mbox{pp. 196--202.} \mbox{(ISBN 980-6560-11-6)}
%\item[{[Ri05]}] Riihentaus, J. ``An integrability condition and weighted boundary behavior of subharmonic and ${\mathcal{M}}$-subharmonic functions:
%a survey'', Int. J. Diff. Eq. Appl., {\textbf{10}} (2005), 1--14.
\item[{[Ri06$_1$]}] Riihentaus, J.  ``A weighted boundary limit result for subharmonic functions'', Adv. Algebra and Analysis, {\textbf{1}}, No.~1 (2006), \mbox{27--38.}
\item[{[Ri06$_2$]}] Riihentaus, J. ``Separately quasi-nearly subharmonic functions'', in:  Complex Analysis and Potential Theory, Conference
Satellite to ICM~2006, Proceedings, Tahir Aliyev Azero$\breve{\textrm{g}}$lu, Promarz M. Tamrazov (eds.), Gebze Institute of Technology, Gebze, Turkey, 
September 8-14,  2006, World Scientific, Singapore, 2007, pp.~156-165.
\item[{[Ri07$_1$]}] Riihentaus, J. ``On the subharmonicity of separately  subharmonic functions'', in: 11th WSEAS International 
Conference on Applied Mathematics (MATH'07),  Dallas, Texas, USA, March 22-24, 2007, Proceedings, Kleanthis Psarris, Andrew D.~Jones (eds.), WSEAS, 2007,  pp.~230-236. \mbox{(IBSN 978-960-8457-60-7)}
\item[{[Ri07$_2$]}] Riihentaus, J. ``On separately harmonic and subharmonic functions'', Int. J. Pure Appl. Math., {\textbf{35}}, No. 4 (2007), 
\mbox{435-446}.
\item[{[Ri07$_3$]}] Riihentaus, J. ``Subharmonic functions, generalizations and separately subharmonic functions'', in:
The XIV-th Conference on Analytic Functions, Che\l m, Poland, July 22-28, 2007,  arXiv:math/0610259v5 [math.AP] 8 Oct 2008, and  Scientific Bulletin of Che\l m, Section of Mathematics and Computer Science, {\textbf{2}} (2007), 49--76. \mbox{(ISBN 978-83-61149-24-8)} 
\item[{[Ri08]}] Riihentaus, J. ``Subharmonic functions, generalizations, weighted boundary behavior,  and separately subharmonic functions: a survey'', 
in: The Fifth World Congress of Nonlinear Analysts (WCNA~2008), Orlando, Florida, USA, July~2-9, 2008, submitted. 
%\item[{[Ru50]}] Rudin, W. ``Integral representation of continuous functions'', Trans. Amer. Math. Soc., {\textbf{68}} (1950), 
%278--286.
%\item[{[Ru79]}] Rudin, W. ``Real and Complex Analysis'', 
%Tata McGraw-Hill, New Delhi, 1979.
%\item[{[Sa41]}] Saks, S. ``On the operators of Blaschke and Privaloff for subharmonic functions'', Rec. Math. (Mat. Sbornik), {\textbf{9 (51)}} (1941), 
%451--456.
%\item[{[Sh56]}] Shapiro, V.L.  ``Generalized laplacians'',  Amer. J. Math.,
%{\textbf{78}} (1956), \mbox{497--508.}
%\item[{[Sh71]}] Shapiro, V.L.  ``Removable sets for pointwise subharmonic functions'', Trans. Amer. Math. Soc.,
%{\textbf{159}} (1971), \mbox{369--380.}
%\item[{[Sh78]}] Shapiro, V.L.  ``Subharmonic functions and Hausdorff measure'', J. Diff. Eq., {\textbf{27}} (1978), \mbox{28--45.}
%\item[{[Si69]}] Siciak, J.  ``Separately analytic functions and envelopes of holomorphy of some lower dimensional subsets of ${\mathbb{C}}^n$'', 
%Ann. Polon. Math.,
%{\textbf{22}} (1969), \mbox{145--171.}
%\item[{[Sz33]}] Szpilrajn, E.  ``Remarques sur les fonctions sousharmoniques'',  Ann. Math., {\textbf{34}} (1933), \mbox{588--594.}
%\item[{[V\" a71]}] V\"ais\"al\"a, J. ``Lectures on n-Dimensional Quasiconformal Mappings'', 
%Lecture Notes in Mathematics 229, Springer-Verlag, Berlin, 1971.  
\item[{[To86]}] Torchinsky, A. ``Real-Variable Methods in Harmonic Analysis'', Academic Press, London, 1986.
\item[{[Vu82]}] Vuorinen, M.  ``On the Harnack constant and the boundary behavior of Harnack functions'',  Ann. Acad. Sci. Fenn., Ser. A I, 
Math., {\textbf{7}}, No.~2 (1982), 259--277.
\item[{[Wi88]}] Wiegerinck, J.  ``Separately  subharmonic functions need not be subharmonic'',  
Proc.  Amer. Math.  Soc., {\textbf{104}}, No.~3 (1988), 770-771.
\item[{[WZ91]}] Wiegerinck, J., Zeinstra, R.  ``Separately subharmonic functions: when are they subharmonic'', in: Proceedings of Symposia in Pure 
Mathematics, Several Complex Variables and Complex Geometry, vol. {\textbf{52}}, part {\textbf{1}},  Eric Bedford, John P. D'Angelo, Robert E.~Greene, Steven G.~Krantz (eds.),  Amer. Math. Soc., Providence, Rhode Island, 1991,  \mbox{pp. 245--249}.
\end{enumerate}
\end{flushleft}
\end{document}